\documentclass[a4paper]{amsart}

\usepackage[mathscr]{eucal}
\usepackage{amsmath}
\usepackage{amsfonts}
\usepackage{amssymb}

\theoremstyle{plain}
\newtheorem{theorem}{Theorem}
\newtheorem{proposition}{Proposition}[section]
\newtheorem{lemma}{Lemma}[section]
\newtheorem{remark}{Remark}

\begin{document}
\title{On the zeros of Eisenstein series for $\Gamma_0^{*} (2)$ and $\Gamma_0^{*} (3)$}
\author{Tsuyoshi Miezaki, Hiroshi Nozaki, Junichi Shigezumi}

\maketitle \vspace{-0.2in}
\begin{center}
Graduate School of Mathematics Kyushu University\\
Hakozaki 6-10-1 Higashi- ku, Fukuoka, 812-8581 Japan\\ \quad
\end{center} \vspace{0.1in}

\begin{quote}
{\small\bfseries Abstract.}
We locate all of the zeros of the Eisenstein series associated with the Fricke groups $\Gamma_0^{*}(2)$ and $\Gamma_0^{*}(3)$ in their fundamental domains by applying and expanding the method of F. K. C. Rankin and H. P. F. Swinnerton-Dyer (``{\it On the zeros of Eisenstein series}'', 1970).\\ \vspace{-0.15in}

\vspace{-0.15in}

\quad
\end{quote}

\section{Introduction}

Let $k \geqslant 4$ be an even integer. For $z \in \mathbb{H} := \{z \in \mathbb{C} \: ; \: Im(z)>0 \}$, let
\begin{equation}
E_k(z) := \frac{1}{2} \sum_{(c,d)=1}(c z + d)^{- k} \label{def:e}
\end{equation}
be the {\it Eisenstein series} associated with $\text{SL}_2(\mathbb{Z})$. Moreover, let
\begin{equation*}
\mathbb{F} := \left\{|z| \geqslant 1, \: - 1 / 2 \leqslant Re(z) \leqslant 0\right\}
 \cup \left\{|z| > 1, \: 0 \leqslant Re(z) < 1 / 2 \right\}
\end{equation*}
be the {\it standaard fundamental domain} for $\text{SL}_2(\mathbb{Z})$.

F. K. C. Rankin and H. P. F. Swinnerton-Dyer considered the problem of locating the zeros of $E_k(z)$ in $\mathbb{F}$ \cite{RSD}. They proved that $n$ zeros are on the arc $A := \{z \in \mathbb{C} \: ; \: |z|=1, \: \pi / 2 < Arg(z) < 2 \pi / 3\}$ for $k = 12 n + s \: (s = 4, 6, 8, 10, 0, \text{ and } 14)$. They also said in the last part of the paper, ``This method can equally well be applied to Eisenstein series associated with subgroup of the modular group.'' However, it seems unclear how widely this claim holds. 

Here, we consider the same problem for Fricke groups $\Gamma_0^{*}(2)$ and $\Gamma_0^{*}(3)$ (See \cite{K}, \cite{Q}), which are commensurable with $\text{SL}_2(\mathbb{Z})$. For a fixed prime $p$, we define the following:
\begin{equation}
\Gamma_0^{*}(p) := \Gamma_0(p) \cup \Gamma_0(p) \: W_p,
\end{equation}
where
\begin{equation}
\Gamma_0(p) := \left\{ \left(\begin{smallmatrix} a & b \\ c & d \end{smallmatrix}\right) \in \text{SL}_2(\mathbb{Z}) \: ; \: c \equiv 0\pmod{p}\right\}, \quad W_p := \left(\begin{smallmatrix} 0 & -1 / \sqrt{p} \\ \sqrt{p} & 0 \end{smallmatrix}\right).
\end{equation}

Let
\begin{equation}
E_{k, p}^{*}(z) := \frac{1}{p^{k / 2}+1} \left(p^{k / 2} E_k(p z) + E_k(z) \right) \label{def:e*}
\end{equation}
be the Eisenstein series associated with $\Gamma_0^{*}(p)$. The regions 
\begin{align*}
\mathbb{F}^{*}(2) &:= \left\{|z| \geqslant 1 / \sqrt{2}, \: - 1 / 2 \leqslant Re(z) \leqslant 0\right\}
 \bigcup \left\{|z| > 1 / \sqrt{2}, \: 0 \leqslant Re(z) < 1 / 2 \right\},\\
\mathbb{F}^{*}(3) &:= \left\{|z| \geqslant 1 / \sqrt{3}, \: - 1 / 2 \leqslant Re(z) \leqslant 0\right\}
 \bigcup \left\{|z| > 1 / \sqrt{3}, \: 0 \leqslant Re(z) < 1 / 2 \right\}
\end{align*}
are fundamental domains for $\Gamma_0^{*}(2)$ and $\Gamma_0^{*}(3)$, respectively.

Define $A_2^{*} := \{z \in \mathbb{C} \: ; \: |z| = 1 / \sqrt{2}, \: \pi / 2 < Arg(z) < 3 \pi / 4\}$, and $A_3^{*} := \{z \in \mathbb{C} \: ; \: |z| = 1 / \sqrt{3}, \: \pi / 2 < Arg(z) < 5 \pi / 6\}$. We then have $\overline{A_2^{*}} = A_2^{*} \cup \{ i / \sqrt{2}, e^{i (3 \pi / 4)} / \sqrt{2} \}$ and $\overline{A_3^{*}} = A_3^{*} \cup \{ i / \sqrt{3}, e^{i (5 \pi / 6)} / \sqrt{3} \}$.

In the present paper, we will apply the method of F. K. C. Rankin and H. P. F. Swinnerton-Dyer (RSD Method) to the Eisenstein series associated with $\Gamma_0^{*}(2)$ and $\Gamma_0^{*}(3)$. We will prove the following theorems:

\begin{theorem}
Let $k \geqslant 4$ be an even integer. All of the zeros of $E_{k, 2}^{*}(z)$ in $\mathbb{F}^{*}(2)$ are on the arc $\overline{A_2^{*}}$. \label{th-g0s2}
\end{theorem}

\begin{theorem}
Let $k \geqslant 4$ be an even integer. All of the zeros of $E_{k, 3}^{*}(z)$ in $\mathbb{F}^{*}(3)$ are on the arc $\overline{A_3^{*}}$. \label{th-g0s3}
\end{theorem}

\section{RSD Method}\label{sect-1}

At the beginning of the proof in \cite{RSD}, F. K. C. Rankin and H. P. F. Swinnerton-Dyer considered the following:
\begin{equation}
F_k(\theta) := e^{i k \theta / 2} E_k\left(e^{i \theta}\right), \label{def:f}
\end{equation}
which is real for all $\theta \in [0, \pi]$. Considering the four terms with $c^2 + d^2 = 1$, they proved that
\begin{equation}
F_k(\theta) = 2 \cos(k \theta / 2) + R_1, \label{eqn-fkt}
\end{equation}
where $R_1$ is the rest of the series ({\it i.e.} $c^2 + d^2 > 1$). Moreover they showed
\begin{equation}
|R_1|
 \leqslant 1 + \left(\frac{1}{2}\right)^{k / 2} + 4 \left(\frac{2}{5}\right)^{k / 2} + \frac{20 \sqrt{2}}{k - 3} \left(\frac{9}{2}\right)^{(3 - k) / 2}. \label{r1bound}
\end{equation}
They computed the value of the right-hand side of (\ref{r1bound}) at $k = 12$ to be approximately $1.03562$, which is monotonically decreasing in $k$. Thus, they could show that $|R_1| < 2$ for all $k \geqslant 12$. If  $\cos (k \theta / 2)$ is $+1$ or $-1$, then $F_k(2 m \pi / k)$ is positive or negative, respectively.

In order to determine the location of all of the zeros of $E_k(z)$ in $\mathbb{F}$, we need the {\it valence formula}:

\begin{proposition}[valence formula]
Let $f$ be a modular function of weight $k$ for $\text{\upshape SL}_2(\mathbb{Z})$, which is not identically zero. We have
\begin{equation}
v_{\infty}(f) + \frac{1}{2} v_{i}(f) + \frac{1}{3} v_{\rho} (f) + \sum_{\begin{subarray}{c} p \in \text{\upshape SL}_2(\mathbb{Z}) \setminus \mathbb{H} \\ p \ne i, \; \rho\end{subarray}} v_p(f) = \frac{k}{12},
\end{equation}
where $v_p(f)$ is the order of $f$ at $p$, and $\rho := e^{i (2 \pi / 3)}$ $($See \cite{S}$)$. \label{prop-vf}
\end{proposition}

Write $m(k) := \left\lfloor \frac{k}{12} - \frac{t}{4} \right\rfloor$, where $t = 0 \text{ or } 2$, such that $t \equiv k \pmod{4}$. Then, $k = 12 m(k) + s \: (s = 4, 6, 8, 10, 0, \text{ and } 14)$.

As F. K. C. Rankin and H. P. F. Swinnerton-Dyer observed, the fact that $E_k(z)$ has $m(k)$ zeros on the arc $A$, the valence formula, and Remark \ref{prop-bd_ord_1} below, imply that all of the zeros of $E_k(z)$ in the standard fundamental domain for $\text{SL}_2(\mathbb{Z})$ are on $A \cup \{ i, \rho \}$ for every even integer $k \geqslant 4$.

\begin{remark}
Let $k \geqslant 4$ be an even integer. We have
\begin{center}
\begin{tabular}{rcccrcc}
$k \pmod{12}$ & $v_{i / \sqrt{3}}(E_k)$ & $v_{\rho_3}(E_k)$ & \qquad & $k \pmod{12}$ & $v_{i / \sqrt{3}}(E_k)$ & $v_{\rho_3}(E_k)$\\
\hline
$0$ & $0$ & $0$ && $6$ & $1$ & $0$\\
$2$ & $1$ & $2$ && $8$ & $0$ & $2$\\
$4$ & $0$ & $1$ && $10$ & $1$ & $1$\\
\hline
\end{tabular}
\end{center}\label{prop-bd_ord_1}
\end{remark}

\section{$\Gamma_0^{*}(2)$ (Proof of Theorem \ref{th-g0s2})}\label{sect-2}

\subsection{Preliminaries}
We define
\begin{equation}
F_{k, 2}^{*}(\theta) := e^{i k \theta / 2} E_{k, 2}^{*}\left(e^{i \theta} / \sqrt{2}\right). \label{def:f*2}
\end{equation}

Before proving Theorem \ref{th-g0s2}, we consider an expansion of $F_{k, 2}^{*}(\theta)$.

By the definition of $E_k(z), E_{k, 2}^{*}(z)$ ({\it cf.} (\ref{def:e}), (\ref{def:e*})), we have
\begin{multline*}
2 (2^{k / 2}+1) e^{i k \theta / 2} E_{k, 2}^{*}\left(e^{i \theta} / \sqrt{2}\right)\\
\quad = 2^{k / 2} \sum_{(c,d)=1}(c e^{-i \theta / 2} + \sqrt{2} d e^{i \theta / 2})^{- k} 
+ 2^{k / 2} \sum_{(c,d)=1}(c e^{i \theta / 2} + \sqrt{2} d e^{-i \theta / 2})^{- k}.
\end{multline*}

Now, $(c,d) = 1$ is split in two cases, namely, $c$ is {\it odd} or $c$ is {\it even}. We consider the case in which $c$ is {\it even}. We have
\begin{align*}
2^{k / 2} \sum_{\begin{subarray}{c} (c,d)=1\\ c:even\end{subarray}}(c e^{-i \theta / 2} + \sqrt{2} d e^{i \theta / 2})^{- k}
&= \sum_{\begin{subarray}{c} (c,d)=1\\ d:odd\end{subarray}}(\sqrt{2} c' e^{-i \theta / 2} + d e^{i \theta / 2})^{- k} \quad (c = 2 c')\\
&= \sum_{\begin{subarray}{c} (c,d)=1\\ c:odd\end{subarray}}(c e^{i \theta / 2} + \sqrt{2} d e^{-i \theta / 2})^{- k}.
\end{align*}

Similarly,
\begin{equation*}
2^{k / 2} \sum_{\begin{subarray}{c} (c,d)=1\\ c:even\end{subarray}}(c e^{i \theta / 2} + \sqrt{2} d e^{-i \theta / 2})^{- k}
 = \sum_{\begin{subarray}{c} (c,d)=1\\ c:odd\end{subarray}}(c e^{-i \theta / 2} + \sqrt{2} d e^{i \theta / 2})^{- k}.
\end{equation*}

Thus, we can write the following:
\begin{equation}
F_{k, 2}^{*}(\theta)
 = \frac{1}{2} \sum_{\begin{subarray}{c} (c,d)=1\\ c:odd\end{subarray}}(c e^{i \theta / 2} + \sqrt{2} d e^{-i \theta / 2})^{- k}
 + \frac{1}{2} \sum_{\begin{subarray}{c} (c,d)=1\\ c:odd\end{subarray}}(c e^{-i \theta / 2} + \sqrt{2} d e^{i \theta / 2})^{- k}.
\end{equation}

Hence, we use this expression as a definition. 

In the last part of this section, we compare the two series in this expression. Note that for any pair $(c,d)$, $(c e^{i \theta / 2} + \sqrt{2} d e^{-i \theta / 2})^{- k}$ and $(c e^{-i \theta / 2} + \sqrt{2} d e^{i \theta / 2})^{- k}$ are conjugates of each other. Thus, we have the following lemma:

\begin{lemma}
$F_{k, 2}^{*}(\theta)$ is real, for all $\theta \in [0, \pi]$. \label{lemma-real2}
\end{lemma}

\subsection{Application of the RSD Method}
We will apply the method of F. K. C. Rankin and H. P. F. Swinnerton-Dyer (RSD Method) to the Eisenstein series associated with $\Gamma_0^{*}(2)$. Note that $N := c^2 + d^2$.

First, we consider the case of $N = 1$. Because $c$ is odd, there are two cases, $(c,d)=(1,0)$ and $(c,d)=(-1,0)$. Then, we can write:
\begin{equation}
F_{k, 2}^{*}(\theta) = 2 \cos(k \theta /2) + R_2^{*},
\end{equation}
where $R_2^{*}$ denotes the remaining terms of the series.

Now, 
\begin{equation*}
|R_2^{*}| \leqslant \sum_{\begin{subarray}{c} (c,d)=1\\ c:odd, \; N > 1\end{subarray}}|c e^{i \theta / 2} + \sqrt{2} d e^{-i \theta / 2}|^{- k}.
\end{equation*}

Let $v_{k}(c,d,\theta) := |c e^{i \theta / 2} + \sqrt{2} d e^{-i \theta / 2}|^{- k}$, then $v_{k}(c,d,\theta) = 1 / \left( c^2 + 2 d^2 + 2 \sqrt{2} c d \cos\theta \right)^{k / 2}$, and $v_{k}(c,d,\theta)=v_{k}(-c,-d,\theta)$.

Next, we will consider the following three cases, namely, $N = 2, 5$, and $N \geqslant 10$. Considering $\theta \in [\pi / 2, 3 \pi / 4]$, we have the following:
\begin{align*}
&\text{When $N = 2$,}&
v_k(1, 1, \theta) &\leqslant 1, \quad
&v_k(1, - 1, \theta) &\leqslant (1 / 3)^{k / 2}.\\
&\text{When $N = 5$,}&
v_k(1, 2, \theta) &\leqslant (1 / 5)^{k / 2},
&v_k(1, - 2, \theta) &\leqslant (1 / 3)^{k}.\\
&\text{When $N \geqslant 10$,}
\end{align*}
\begin{align*}
|c e^{i \theta / 2} \pm \sqrt{2} d e^{-i \theta / 2}|^2
&\geqslant c^2 + 2 d^2 - 2 \sqrt{2} |c d| |\cos\theta|\\
&\geqslant (c^2 + d^2) / N \quad = N / 3,
\end{align*}
and the remaining problem concerns the number of terms with $c^2 + d^2 = N$. Because $c$ is odd, $|c| = 1, 3, ... , 2 N' - 1 \leqslant N^{1/2}$, so the number of $|c|$ is not more than $(N^{1/2}+1)/ 2$. Thus, the number of terms with $c^2 + d^2 = N$ is not more than $2 (N^{1/2}+1) \leqslant 3 N^{1/2}$, for $N \geqslant 5$. Then,
\begin{align*}
\sum_{\begin{subarray}{c} (c,d)=1\\ c:odd, \; N \geq 10\end{subarray}}|c e^{i \theta / 2} + \sqrt{2} d e^{-i \theta / 2}|^{- k}
&\leqslant \sum_{N=10}^{\infty} 3 N^{1/2} \left(\frac{N}{3}\right)^{- k / 2}\\
 &\leqslant \frac{18 \sqrt{3}}{k-3} \left(\frac{1}{3}\right)^{(k-3)/2}
 = \frac{162}{k-3} \left(\frac{1}{3}\right)^{k / 2}.
\end{align*}
Thus,
\begin{equation}
|R_2^{*}|
 \leqslant 2 + 2 \left(\frac{1}{3}\right)^{k / 2} + 2 \left(\frac{1}{5}\right)^{k / 2} + 2 \left(\frac{1}{3}\right)^{k} + \frac{162}{k-3} \left(\frac{1}{3}\right)^{k / 2}. \label{r*2bound0}
\end{equation}

Recalling the previous section (RSD Method), we want to show that $|R_2^{*}| < 2$. However, the right-hand side is greater than 2, so this bound is not good. The case in which $(c,d) = \pm (1,1)$ gives a bound equal to 2. We will consider the expansion of the method in the following sections.

\subsection{Expansion of the RSD Method (1)}
In the previous subsection, we could not obtain a good bound for $|R_2^{*}|$, where $(c,d) = \pm (1,1)$. Note that ``$v_k(1,1,\theta) = 1 \Leftrightarrow \theta = 3 \pi / 4$''. Furthermore, ``$v_k(1,1,\theta) < 1 \Leftrightarrow \theta < 3 \pi / 4$''. Therefore, we can easily expect that a good bound can be obtained for $\theta \in [\pi / 2, 3 \pi / 4 - x]$ for small $x > 0$. However, if $k = 8 n$, then we need $|R_2^{*}| < 2$ for $\theta = 3 \pi / 4$ in this method. In the next section, we will consider the case in which $k = 8 n, \theta = 3 \pi / 4$.

Define $m_2(k) := \left\lfloor \frac{k}{8} - \frac{t}{4} \right\rfloor$, where $t=0, 2$ is chosen so that $t \equiv k \pmod{4}$, and $\lfloor n \rfloor$ is the largest integer not more than $n$.

Let $k = 8 n + s \: (n = m_2(k), \: s = 4, 6, 0, \text{ and } 10)$. We may assume that $k \geqslant 8$.

The first step is to consider how small $x$ should be. We consider each of the cases $s = 4, 6, 0, \: and \: 10$.

When $s = 4$, for $\pi / 2 \leqslant \theta \leqslant 3 \pi / 4$, $(2 n + 1) \pi \leqslant k \theta / 2 \; (= (4 n + 2) \theta) \leqslant (3 n + 1) \pi + \pi / 2$. So the last integer point ({\it i.e.} $\pm 1$) is $k \theta / 2 = (3 n + 1) \pi$, then $\theta = 3 \pi / 4 - \pi / k$. Similarly, when $s = 6, \: and \; 10$, the last integer points are $\theta = 3 \pi / 4 - \pi / 2 k, \: 3 \pi / 4 - 3 \pi / 2 k$, respectively. When $s = 0$, the second to the last integer point is $\theta = 3 \pi / 4 - \pi / k$.

Thus, we need $x \leqslant \pi / 2 k$.

\begin{lemma}
Let $k \geqslant 8$. For all $\theta \in [\pi / 2, 3 \pi / 4 - x] \; (x = \pi / 2 k)$, $|R_2^{*}| < 2$. \label{lemma-r*2}
\end{lemma}

\begin{proof}
Let $k \geqslant 8$ and $x = \pi / 2 k$, then $0 \leqslant x \leqslant \pi / 16$. If $0 \leqslant x \leqslant \pi / 16$, then $1 - \cos x \geqslant \frac{31}{64} x^2$.
\begin{align*}
|e^{i \theta / 2} + \sqrt{2} e^{-i \theta / 2}|^2
&\geqslant 3 + 2 \sqrt{2} \cos(3 \pi / 4 - x)= 1 + 2 (1 - \cos x) + 2 \sin x\\
 &\geqslant 1 + 4 (1 - \cos x) \geqslant 1 + (31 / 16) x^2.
\end{align*}
\begin{equation*}
|e^{i \theta / 2} + \sqrt{2} e^{-i \theta / 2}|^k \geqslant \left(1 + (31 / 16) x^2\right)^{k / 2} \geqslant 1 + (31 / 4) x^2. \quad (k \geqslant 8)
\end{equation*}
\begin{equation*}
v_k(1, 1, \theta)
 \leqslant \frac{1}{1 + (31 / 4) x^2}
 \leqslant 1 - \frac{31 \times 256}{31 \pi^2 + 1024} x^2.
\end{equation*}
Thus,
\begin{equation*}
2 v_k(1, 1, \theta) \leqslant 2 - (265 / 9) / k^2.
\end{equation*}

Furthermore,
\begin{equation*}
2 \left(\frac{1}{3}\right)^{k / 2} + 2 \left(\frac{1}{5}\right)^{k / 2} + 2 \left(\frac{1}{3}\right)^{k} + \frac{162}{k-3} \left(\frac{1}{3}\right)^{k / 2} \leqslant 35 \left(\frac{1}{3}\right)^{k / 2} \quad (k \geqslant 8).
\end{equation*}

Then, we have
\begin{equation*}
|R_2^{*}| \leqslant 2 - \frac{265}{9} \frac{1}{k^2} + 35 \left(\frac{1}{3}\right)^{k}.
\end{equation*}

Next, if we can show that
\begin{equation*}
35 \left(\frac{1}{3}\right)^{k / 2} < \frac{265}{9} \frac{1}{k^2} \quad or \quad \frac{3^{k / 2}}{35} > \frac{9}{265} k^2,
\end{equation*}
then the bound is less than $2$. The proof will thus be complete.

Let $f(x) := (1 / 35) 3^{x/2} - \frac{9}{265} x^2$. Then, $f'(x) = (\log 3 / 70) 3^{x/2} - \frac{18}{265} x$, $f''(x) = ((\log 3)^2 / 140) 3^{x/2} - \frac{18}{265}$. First, $f''$ is monotonically increasing for  $x \geqslant 8$, and $f''(8) = 0.63038... > 0$, so $f'' > 0$ for  $x \geqslant 8$. Second, $f'$ is monotonically increasing for  $x \geqslant 8$, and $f'(8) = 0.72785... > 0$, so $f' > 0$ for  $x \geqslant 8$. Finally, $f$ is monotonically increasing for $x \geqslant 8$, and $f(8) = 0.14070... > 0$, so $f > 0$ for  $x \geqslant 8$.
\end{proof}

\subsection{Expansion of the RSD Method (2)}
For the case of ``$k = 8 n, \theta = 3 \pi / 4$'', we need the following lemma:
\begin{lemma}
Let $k$ be an integer such that $k = 8 n$ for some $n \in \mathbb{N}$. If $n$ is even, then $F_{k, 2}^{*}(3 \pi / 4) > 0$. On the other hand if $n$ is odd, then $F_{k, 2}^{*}(3 \pi / 4) < 0$. \label{lemma-8n}
\end{lemma}

\begin{proof}
Let $k = 8 n$ ($n \geqslant 1$). By the definition of $E_{k, 2}^{*}(z), F_{k, 2}^{*}(z)$ ({\it cf.} (\ref{def:e*}), (\ref{def:f*2})), we have
\begin{equation*}
F_{k, 2}^{*}(3 \pi / 4) = \frac{e^{i 3 (k / 8) \pi}}{2^{k / 2}+1} \left( 2^{k / 2} E_k(- 1 + i) + E_k\left(\frac{- 1 + i}{2}\right) \right).
\end{equation*}

By the {\it transformation rule} for $\text{SL}_2(\mathbb{Z})$,
\begin{equation*}
E_k(- 1 + i) = E_k(i), \quad
E_k\left( (- 1 + i) / 2 \right) = (1 + i)^k E_k(1 + i) = 2^{k / 2} E_k(i).
\end{equation*}

Then,
\begin{equation}
F_{8n,2}^{*}(3 \pi / 4)
 = 2 e^{i n \pi} \frac{2^{4n}}{2^{4n} + 1} F_{8n}(\pi / 2), \label{eq-2ex2}
\end{equation}
where $\frac{2^{4n}}{2^{4n} + 1} > 0$, $F_{8n}(\pi / 2) = 2 \cos(2 n \pi) + R_1 > 0$. The question is then: ``Which holds, $F_k(\pi / 2) < 0$ or $F_k(\pi / 2) > 0$?''

F. K. C. Rankin and H. P. F. Swinnerton-Dyer showed (\ref{eqn-fkt}) and (\ref{r1bound}) \cite{RSD}. They then proved that $|R_1| < 2$ for $k \geqslant 12$. This was necessary only for $k \geqslant 12$. Now we need $|R_1| < 2$ for $k \geqslant 8$. The value of the right-hand side of (\ref{r1bound}) at $k = 8$ is $1.29658... < 2$, which is monotonically decreasing in $k$. Thus, we can show 
\begin{equation}
|R_1| < 2 \qquad \text{for all} \quad k \geqslant 8. \label{r1newbound}
\end{equation}
Then, the sign ($\pm$) of $F_{k, 2}^{*}(3 \pi / 4)$ is that of $e^{i n \pi}$. Thus, the proof is complete.
\end{proof}

Next, we proved that $E_{k, 2}^{*}(z)$ has $m_2(k)$ zeros on the arc $A_2^{*}$. In order to determine the location of all of the zeros of $E_{k, 2}^{*}(z)$ in $\mathbb{F}^{*}(2)$, we need the valence formula for $\Gamma_0^{*}(2)$:

\begin{proposition}
Let $f$ be a modular function of weight $k$ for $\Gamma_0^{*}(2)$, which is not identically zero. We have
\begin{equation}
v_{\infty}(f) + \frac{1}{2} v_{i / \sqrt{2}}(f) + \frac{1}{4} v_{\rho_2} (f) + \sum_{\begin{subarray}{c} p \in \Gamma_0^{*}(2) \setminus \mathbb{H} \\ p \ne i / \sqrt{2}, \; \rho_2\end{subarray}} v_p(f) = \frac{k}{8},
\end{equation}
where $\rho_2 := e^{i (3 \pi / 4)} \big/ \sqrt{2}$. \label{prop-vf-g0s2}
\end{proposition}

The proof of this proposition is similar to that for Proposition \ref{prop-vf} (See \cite{S}).

If $k \equiv 4,6, \text{ and } 0 \pmod{8}$, then $k / 8 - m_2(k) < 1$. Thus, all of the zeros of $E_{k, 2}^{*}(z)$ in $\mathbb{F}^{*}(2)$ are on the arc $\overline{A_2^{*}}$. On the other hand, if $k \equiv 2 \pmod{8}$, then we have $E_{k, 2}^{*}(i / \sqrt{2}) = i^k E_{k, 2}^{*}(i / \sqrt{2})$ by the transformation rule for $\Gamma_0^{*}(2)$. Then, we have $k / 8 - m_2(k) - v_{i / \sqrt{2}}(E_{k, 2}^{*}) / 2 < 1$.

In conclusion, for every even integer $k \geqslant 4$, all of the zeros of $E_{k, 2}^{*}(z)$ in $\mathbb{F}^{*}(2)$ are on the arc $\overline{A_2^{*}}$.

\begin{remark}$($See Proposition \ref{prop-bd_ord_2}$)$
Let $k \geqslant 4$ be an even integer. We have
\begin{center}
\begin{tabular}{rcccrcc}
$k \pmod{8}$ & $v_{i / \sqrt{2}}(E_{k, 2}^{*})$ & $v_{\rho_2}(E_{k, 2}^{*})$ & \qquad & $k \pmod{8}$ & $v_{i / \sqrt{2}}(E_{k, 2}^{*})$ & $v_{\rho_2}(E_{k, 2}^{*})$\\
\hline
$0$ & $0$ & $0$ && $4$ & $0$ & $2$\\
$2$ & $1$ & $3$ && $6$ & $1$ & $1$\\
\hline
\end{tabular}
\end{center}
\end{remark}\quad

\section{$\Gamma_0^{*}(3)$ (Proof of Theorem \ref{th-g0s3})}\label{sect-3}

\subsection{Preliminaries}
We define
\begin{equation}
F_{k, 3}^{*}(\theta) := e^{i k \theta / 2} E_{k, 3}^{*}\left(e^{i \theta} / \sqrt{3}\right). \label{def:f*3}
\end{equation}

Similar to the case of $F_{k, 2}^{*}(\theta)$, we can write the following:
\begin{equation}
F_{k, 3}^{*}(\theta)
 = \frac{1}{2} \sum_{\begin{subarray}{c} (c,d)=1\\ 3 \nmid c\end{subarray}}(c e^{i \theta / 2} + \sqrt{3} d e^{-i \theta / 2})^{- k}
 + \frac{1}{2} \sum_{\begin{subarray}{c} (c,d)=1\\ 3 \nmid c\end{subarray}}(c e^{-i \theta / 2} + \sqrt{3} d e^{i \theta / 2})^{- k}.
\end{equation}

The following lemma is then obtained:

\begin{lemma}
$F_{k, 3}^{*}(\theta)$ is real, for all $\theta \in [0, \pi]$. \label{lemma-real3}
\end{lemma}

\subsection{Application of the RSD Method}
Note that $N := c^2 + d^2$.

First, we consider the case of $N = 1$. We can then write the following:
\begin{equation}
F_{k, 3}^{*}(\theta) = 2 \cos(k \theta /2) + R_3^{*},
\end{equation}
where $R_3^{*}$ denotes the remaining terms.

Let $v_{k}(c, d, \theta) := |c e^{i \theta / 2} + \sqrt{3} d e^{-i \theta / 2}|^{- k}$. We will consider the following cases: $N = 2, 5, 10, 13, 17$, and $N \geqslant 25$. Considering $\theta \in [\pi / 2, 5 \pi / 6]$, we have the following:
\begin{allowdisplaybreaks}
\begin{align*}
&\text{When $N = 2$,}&
v_k(1, 1, \theta) &\leqslant 1, \quad
&v_k(1, - 1, \theta) &\leqslant (1 / 2)^{k}.\\
&\text{When $N = 5$,}&
v_k(1, 2, \theta) &\leqslant (1 / 7)^{k / 2},
&v_k(1, - 2, \theta) &\leqslant (1 / 13)^{k / 2},\\
&&v_k(2, 1, \theta) &\leqslant 1,
&v_k(2, - 1, \theta) &\leqslant (1 / 7)^{k / 2}.\\
&\text{When $N = 10$,}&
v_k(1, 3, \theta) &\leqslant (1 / 19)^{k / 2},
&v_k(1, - 3, \theta) &\leqslant (1 / 28)^{k / 2}.\\
&\text{When $N = 13$,}&
v_k(2, 3, \theta) &\leqslant (1 / 13)^{k / 2},
&v_k(2, - 3, \theta) &\leqslant (1 / 31)^{k / 2}.\\
&\text{When $N = 17$,}&
v_k(1, 4, \theta) &\leqslant (1 / 37)^{k / 2},
&v_k(1, - 4, \theta) &\leqslant (1 / 7)^{k},\\
&&v_k(4, 1, \theta) &\leqslant (1 / 7)^{k / 2},
&v_k(4, - 1, \theta) &\leqslant (1 / 19)^{k / 2}.\\
&\text{When $N \geqslant 25$,}&
|c e^{i \theta / 2} \pm \sqrt{3} d e^{-i \theta / 2}|^2 \geqslant& N / 6,
\end{align*}
\end{allowdisplaybreaks}
and the number of terms with $c^2 + d^2 = N$ is at most $(11 / 3) N^{1/2}$, for $N \geqslant 16$. Then,
\begin{equation*}
\sum_{\begin{subarray}{c} (c,d)=1\\ 3 \nmid c, \; N \geq 25\end{subarray}}|c e^{i \theta / 2} + \sqrt{3} d e^{-i \theta / 2}|^{- k}
 \leqslant \sum_{N=25}^{\infty} \frac{11}{3} N^{1/2} \left(\frac{1}{6} N\right)^{- k / 2}
 \leqslant \frac{352 \sqrt{6}}{k-3} \left(\frac{1}{2}\right)^{k}.
\end{equation*}
Thus,
\begin{equation}
|R_3^{*}|
 \leqslant 4
 + 2 \left(\frac{1}{2}\right)^{k} + 6 \left(\frac{1}{7}\right)^{k / 2}
 \cdots + 2 \left(\frac{1}{7}\right)^{k}
 + \frac{352 \sqrt{6}}{k-3} \left(\frac{1}{2}\right)^{k}. \label{r*3bound0}
\end{equation}

The cases of $(c, d) = \pm (1, 1), \: \pm (2, 1)$ give a bound equal to 4. We will consider an expansion of the method similar to that of $\Gamma_0^{*}(2)$.

\subsection{Expansion of the RSD Method (1)}
Similar to the method of $\Gamma_0^{*}(2)$, we will consider $\theta \in [\pi / 2, 5 \pi / 6 - x]$ for small $x > 0$. In the next subsection, we also consider the case in which $k = 12 n, \: \theta = 5 \pi / 6$.

Define $m_3(k) := \left\lfloor \frac{k}{6} - \frac{t}{4} \right\rfloor$, where $t=0, 2$ is chosen so that $t \equiv k \pmod{4}$. We may assume that $k \geqslant 8$.

How small should $x$ be? Let $k = 12 m_3(k) + s$. Considering each case, namely, $s = 4, 6, 8, 10 \: and \; 14$, we need $x \leqslant \pi / 3 k$. 

\begin{lemma}
Let $k \geqslant 8$. For all $\theta \in [\pi / 2, 5 \pi / 6 - x] \; (x = \pi / 3 k)$, $|R_3^{*}| < 2$. \label{lemma-r*3}
\end{lemma}

Before proving the above lemma, we need the following preliminaries.

\begin{proposition}
Let $k \geqslant 8$ be an even integer and $x = \pi / 3 k$, then
\begin{equation*}
4 + 2 \sqrt{3} \cos \left(\frac{5 \pi}{6} - x\right) \geqslant \left(\frac{3}{2}\right)^{2 / k} \left(1 + \frac{256 \times 7 \times 13}{3 \times 127 \times k} x^2\right).
\end{equation*}
\label{prop-cos1}
\end{proposition}

\begin{proof}
We have,
\begin{equation*}
\left(\frac{3}{2}\right)^{2 / k}
 = \sum_{n = 0}^{\infty} \frac{(2 \log 3 / 2)^n}{n!} \frac{1}{k^n}
 \leqslant 1 + \left(2 \log \frac{3}{2}\right) \frac{1}{k} + \frac{1}{2} \left(2 \log \frac{3}{2}\right)^2 \left(\frac{3}{2}\right)^{2 / k} \frac{1}{k^2},
\end{equation*}
\begin{equation*}
3 + 2 \sqrt{3} \cos \left(\frac{5 \pi}{6} - \frac{\pi}{3k}\right) \geqslant \frac{\pi}{\sqrt{3}} \frac{1}{k}.
\end{equation*}

Let $x = \pi / 3 k$. Then, we have
\begin{equation*}
f_1(k) := 4 + 2 \sqrt{3} \cos \left(\frac{5 \pi}{6} - \frac{\pi}{3 k}\right) - \left(\frac{3}{2}\right)^{2 / k} \left(1 + \frac{256 \times 7 \times 13 \times \pi^2}{27 \times 127} \frac{1}{k^3}\right).
\end{equation*}
If $k = 8$, then $f_1(8) = 0.00012876... > 0$. Next, if $k \geqslant 10$, then
\begin{align*}
f_1(k)
 &\geqslant \frac{1}{k} \left\{\frac{\pi}{\sqrt{3}} - 2 \log \frac{3}{2} - \frac{1}{2} \left(2 \log \frac{3}{2}\right)^2 \left(\frac{3}{2}\right)^{2 / k} \frac{1}{k} - \frac{256 \times 7 \times 13 \times \pi^2}{27 \times 127} \left(\frac{3}{2}\right)^{2 / k} \frac{1}{k^2}\right\}\\
 &\geqslant  \frac{1}{k} \times 0.24004... \quad (k \geqslant 10) \quad > 0.
\end{align*}
\end{proof}

\begin{proposition}Let $k \geqslant 8$ be an even integer and $x = \pi / 3 k$, then
\begin{equation*}
7 + 4 \sqrt{3} \cos \left(\frac{5 \pi}{6} - x\right) \geqslant 3^{2 / k} \left(1 + \frac{256 \times 7 \times 13}{3 \times 127 \times k} x^2\right).
\end{equation*}
\label{prop-cos2}
\end{proposition}

\begin{proof}
We have
\begin{gather*}
3^{2 / k} \leqslant 1 + (2 \log 3) \frac{1}{k} + \frac{1}{2} (2 \log 3)^2 3^{2 / k} \frac{1}{k^2},\\
6 + 4 \sqrt{3} \cos \left(\frac{5 \pi}{6} - \frac{\pi}{3k}\right) \geqslant \frac{2 \pi}{\sqrt{3}} \frac{1}{k}.
\end{gather*}

Similar to the proof of Proposition \ref{prop-cos1}, let $x = \pi / 3 k$, and write
\begin{equation*}
f_2(k) := 7 + 4 \sqrt{3} \cos \left(\frac{5 \pi}{6} - \frac{\pi}{3 k}\right) - 3^{2 / k} \left(1 + \frac{256 \times 7 \times 13 \times \pi^2}{27 \times 127} \frac{1}{k^3}\right).
\end{equation*}
If $k = 8$, then $f_2(8) = 0.015057... > 0$. Next, if $k \geqslant 10$, then
\begin{align*}
f_2(k) \geqslant  \frac{1}{k} \times 0.29437... \quad > 0.
\end{align*}
\end{proof}

\begin{proof}[Proof of Lemma \ref{lemma-r*3}]
Let $k \geqslant 8$ and $x = \pi / 3 k$, then $0 \leqslant x \leqslant \pi / 24$.

By Proposition \ref{prop-cos1}
\begin{equation*}
|e^{i \theta / 2} + \sqrt{3} e^{-i \theta / 2}|^2
 \geqslant \left(\frac{3}{2}\right)^{2 / k} \left(1 + \frac{256 \times 7 \times 13}{3 \times 127 \times k} x^2\right).
\end{equation*}
\begin{align*}
|e^{i \theta / 2} + \sqrt{3} e^{-i \theta / 2}|^k
 &\geqslant \left(\frac{3}{2}\right) \left(1 + \frac{256 \times 7 \times 13}{3 \times 127 \times k} x^2\right)^{k / 2}\\
 &\geqslant \frac{3}{2} + \frac{64 \times 7 \times 13}{127} x^2. \quad (k \geqslant 8)
\end{align*}
\begin{align*}
v_k(1, 1, \theta)
&\leqslant \frac{2}{3} - \frac{(128 \times 7 \times 13 / 127)}{(9 / 2) + (64 \times 3 \times 7 \times 13 / 127) x^2} x^2
 \leqslant \frac{2}{3} - \frac{107}{8} x^2. \quad (x \leqslant \pi / 24)
\end{align*}

Similarly, by Proposition \ref{prop-cos2}
\begin{equation*}
|2 e^{i \theta / 2} + \sqrt{3} e^{-i \theta / 2}|^2 \geqslant 3^{2 / k} \left(1 + \frac{256 \times 7 \times 13}{3 \times 127 \times k} x^2\right).
\end{equation*}
\begin{equation*}
|2 e^{i \theta / 2} + \sqrt{3} e^{-i \theta / 2}|^k \geqslant 3 + \frac{128 \times 7 \times 13}{127} x^2.
\end{equation*}
\begin{equation*}
v_k(2, 1, \theta) \leqslant \frac{1}{3} - \frac{107}{16} x^2.
\end{equation*}

Thus,
\begin{equation*}
2 v_k(1, 1, \theta) + 2 v_k(2, 1, \theta) \leqslant 2 - \frac{107 \pi^2}{24} \frac{1}{k^2}.
\end{equation*}

In Eq. (\ref{r*3bound0}), replace $4$ with the bound $2 - \frac{107 \pi^2}{24} \frac{1}{k^2}$. Then,
\begin{equation*}
|R_3^{*}| \leqslant 2 - \frac{107 \pi^2}{24} \frac{1}{k^2}  + 176 \left(\frac{1}{2}\right)^{k}.
\end{equation*}

Similarly to the method for $\Gamma_0^{*}(2)$, we can easily show that the bound is less than two for $k \geqslant 8$.
\end{proof}

\subsection{Expansion of the RSD Method (2)}
For the case ``$k = 12 n, \theta = 5 \pi / 6$'', we need the following lemma:
\begin{lemma}
Let $k$ be the integer such that $k = 12 n$ for some $n \in \mathbb{N}$. If $n$ is even, then $F_{k, 3}^{*}(5 \pi / 6) > 0$. On the other hand, if $n$ is odd, then $F_{k, 3}^{*}(5 \pi / 6) < 0$. \label{lemma-12n}
\end{lemma}

\begin{proof}
Let $k = 12 n$ ($n \geqslant 1$). Similarly to (\ref{eq-2ex2}),
\begin{equation}
F_{12n,3}^{*}(5 \pi / 6)
 = 2 e^{i n \pi} \frac{3^{6n}}{3^{6n} + 1} F_{12n}(2 \pi / 3),
\end{equation}
where $\frac{3^{6n}}{3^{6n} + 1} > 0$, $F_{12n}(2 \pi / 3) = 2 \cos(4 n \pi) + R_1 > 0$ ({\it cf.} (\ref{r1newbound})).
\end{proof}

\begin{proposition}[Valence formula for $\Gamma_0^{*}(3)$]
Let $f$ be a modular function of weight $k$ for $\Gamma_0^{*}(3)$, which is not identically zero. We have
\begin{equation}
v_{\infty}(f) + \frac{1}{2} v_{i / \sqrt{3}}(f) + \frac{1}{6} v_{\rho_3} (f) + \sum_{\begin{subarray}{c} p \in \Gamma_0^{*}(3) \setminus \mathbb{H} \\ p \ne i / \sqrt{3}, \; \rho_3\end{subarray}} v_p(f) = \frac{k}{6},
\end{equation}
where $\rho_3 := e^{i (5 \pi / 6)} \big/ \sqrt{3}$. \label{prop-vf-g0s3}
\end{proposition}

If $k \equiv 4,8,10 \text{ and } 0 \pmod{12}$, then $k / 6 - m_3(k) < 1$. On the other hand, if $k \equiv 2, \: 6 \pmod{12}$, we have $E_{k, 3}^{*}(i / \sqrt{3}) = 0$ and $k / 6 - m_3(k) - v_{i / \sqrt{3}}(E_{k, 3}^{*}) / 2 < 1$.

In conclusion, for every even integer $k \geqslant 4$, all of the zeros of $E_{k, 3}^{*}(z)$ in $\mathbb{F}^{*}(3)$ are on the arc $\overline{A_3^{*}}$.

\begin{remark}$($See Proposition \ref{prop-bd_ord_3}$)$
Let $k \geqslant 4$ be an even integer. We have
\begin{center}
\begin{tabular}{rcccrcc}
$k \pmod{12}$ & $v_{i / \sqrt{3}}(E_{k, 3}^{*})$ & $v_{\rho_3}(E_{k, 3}^{*})$ & \qquad & $k \pmod{12}$ & $v_{i / \sqrt{3}}(E_{k, 3}^{*})$ & $v_{\rho_3}(E_{k, 3}^{*})$\\
\hline
$0$ & $0$ & $0$ && $6$ & $1$ & $3$\\
$2$ & $1$ & $5$ && $8$ & $0$ & $2$\\
$4$ & $0$ & $4$ && $10$ & $1$ & $1$\\
\hline
\end{tabular}
\end{center}
\end{remark}\quad

\begin{remark}
Getz\cite{G} considered a similar problem for the zeros of extremal modular forms of $\text{\upshape SL}_2(\mathbb{Z})$. It seems that similar results do not hold for extremal modular forms of $\Gamma_0^{*}(2)$ and $\Gamma_0^{*}(3)$. We plan to look into this in the near future.
\end{remark}\quad

\appendix

\begin{center}
{\bfseries APPENDIX.}\quad{\bfseries On the space of modular forms for $\Gamma_0^{*} (2)$ and $\Gamma_0^{*} (3)$}
\end{center}

We need theories of the spaces of modular forms for $\Gamma_0^{*} (2)$ and $\Gamma_0^{*} (3)$ in order to decide the orders at some zeros. We refer to {\it J. -P. Serre's ``A Course in Arithmetic''} \cite{S}, which presents theories for the space of modular forms for $\text{SL}_2(\mathbb{Z})$.

Let $M_{k, p}$ be the space of modular forms for $\Gamma_0^{*} (p)$ of weight $k$, and let $M_{k, p}^0$ be the space of cusp forms for $\Gamma_0^{*} (p)$ of weight $k$. When we consider the map $M_{k, p} \ni f \mapsto f(\infty) \in \mathbb{C}$, the kernel of the map is $M_{k, p}^0$. So $\dim(M_{k, p} / M_{k, p}^0) \leqslant 1$, and $M_{k, p} = \mathbb{C} E_{k, p}^{*} \oplus M_{k, p}^0$.

\section{$\Gamma_0^{*} (2)$}

\begin{theorem}Let $k$ be an even integer, and let $\Delta_2 := \frac{17}{1152} ((E_{4,2}^{*})^2-E_{8,2}^{*})$.
\def\labelenumi{(\arabic{enumi})}
\begin{enumerate}
\item
 For $k < 0$ and $k = 2$, $M_{k, 2} = 0$.
\item
 For $k = 0, 4, 6$, and $10$, we have $M_{k, 2}^0 = 0$, and $\dim(M_{k, 2}) = 1$ with a base $E_{k, 2}^{*}$.
\item
 $M_{k, 2}^0 = \Delta_2 M_{k - 8, 2}$.
\end{enumerate}
\def\labelenumi{\arabic{enumi}.}
\end{theorem}

We can prove the above theorem in a similar manner to the proof of to Theorem 4 of \cite[Chapter VII]{S}. We use the valence formula for $\Gamma_0^{*} (2)$ (Proposition \ref{prop-vf-g0s2}).

Furthermore, for a non-negative integer $k$, $\dim(M_{k, 2}) = \lfloor k / 8 \rfloor$ if $k \equiv 2 \pmod{8}$, and $\dim(M_{k, 2}) = \lfloor k / 8 \rfloor + 1$ if $k \not\equiv 2 \pmod{8}$. We have $M_{k, 2} = \mathbb{C} E_{k - 8 n, 2}^{*} E_{8 n, 2}^{*} \oplus M_{k, 2}^0$. Then,
\begin{equation*}
M_{k, 2} = E_{k - 8 n, 2}^{*} (\mathbb{C} E_{8 n, 2}^{*} \oplus \mathbb{C} E_{8 (n - 1), 2}^{*} \Delta_2 \oplus \cdots \oplus \mathbb{C} \Delta_2^n)
\end{equation*}
Thus, for every $p \in \mathbb{H}$ and for every $f \in M_{k, 2}$, $v_p(f) \geqslant v_p(E_{k - 8 n, 2}^{*})$. We also have $E_{10, 2}^{*} = E_{4, 2}^{*} E_{6, 2}^{*}$.

Finally, we have the following proposition:
\begin{proposition}
Let $k \geqslant 4$ be an even integer. For every $f \in M_{k, 2}$, we have
\begin{equation}
\begin{split}
v_{i / \sqrt{2}}(f) \geqslant s_k &\quad(s_k=0, 1 \; \text{such that} \; 2 s_k \equiv k \pmod{4}),\\
v_{\rho_2}(f) \geqslant t_k &\quad(t_k=0, 1, 2, 3 \; \text{such that} \; - 2 t_k \equiv k \pmod{8}).
\end{split}
\end{equation}
In particular, if $f$ is a constant multiple of $E_{k, 2}^{*}$, then the equalities hold. \label{prop-bd_ord_2}
\end{proposition}

\section{$\Gamma_0^{*} (3)$}

\begin{theorem}Let $k$ be an even integer.
\def\labelenumi{(\arabic{enumi})}
\begin{enumerate}
\item
 For $k < 0$ and $k = 2$, $M_{k, 3} = 0$.
\item
 For $k = 0, 4, 6$, we have $M_{k, 3}^0 = 0$, and $\dim(M_{k, 3}) = 1$ with a base $E_{k, 3}^{*}$.
\item
 For $k = 8, 10, 14$, we have $M_{k, 3}^0 = \mathbb{C} \Delta_{3, k}$.
\item
 For $k = 12$, we have $M_{12, 3}^0 = \mathbb{C} \Delta_{3, 12}^0 \oplus \mathbb{C} \Delta_{3, 12}^1$.
\item
 $M_{k, 3}^0 = M_{12, 3}^0 M_{k - 12, 3}$.
\end{enumerate}
\def\labelenumi{\arabic{enumi}.}
where $\Delta_{3, 8} := \frac{41}{1728} ((E_{4, 3}^{*})^2 - E_{8, 3}^{*})$, $\Delta_{3, 10} := \frac{61}{432} (E_{4, 3}^{*} E_{6, 3}^{*} - E_{10, 3}^{*})$, $\Delta_{3, 12}^0 := (\Delta_{3, 8})^2 / E_{4,3}^{*}$, $\Delta_{3, 12}^1 := \Delta_{3, 8} E_{4, 3}^{*}$, and $\Delta_{3, 14} := \Delta_{3, 10} E_{4, 3}^{*}$.
 \label{th-mod_sp_3}
\end{theorem}

Now, we have the following table:

\begin{center}
\begin{tabular}{rccccccrccccc}
$k$ & $f$ & $v_{\infty}$ & $v_{i / \sqrt{3}}$ & $v_{\rho_3}$ & $V_3^{*}$ &\qquad& $k$ & $f$ & $v_{\infty}$ & $v_{i / \sqrt{3}}$ & $v_{\rho_3}$ & $V_3^{*}$\\
\hline
$4$ & $E_{4, 3}^{*}$ & $0$ & $0$ & $4$ & $0$ && $12$ & $E_{12, 3}^{*}$ & $0$ & $0$ & $0$ & $2$\\
$6$ & $E_{6, 3}^{*}$ & $0$ & $1$ & $3$ & $0$ && & $\Delta_{3, 12}^0$ & $2$ & $0$ & $0$ & $0$\\
$8$ & $E_{8, 3}^{*}$ & $0$ & $0$ & $2$ & $1$ && & $\Delta_{3, 12}^1$ & $1$ & $0$ & $6$ & $0$\\
 & $\Delta_{3, 8}$ & $1$ & $0$ & $2$ & $0$ && $14$ & $E_{14, 3}^{*}$ & $0$ & $1$ & $5$ & $1$\\
$10$ & $E_{10, 3}^{*}$ & $0$ & $1$ & $1$ & $1$ && & $\Delta_{3, 14}$ & $1$ & $1$ & $5$ & $0$\\
 & $\Delta_{3, 10}$ & $1$ & $1$ & $1$ & $0$ && & & & & & \\
\hline
\end{tabular}
\end{center}
where $V_3^{*}$ denotes the number of simple zeros of $f$ on $A_3^{*}$.

Furthermore, for a non-negative integer $k$, $\dim(M_{k, 3}) = \lfloor k / 6 \rfloor$ if $k \equiv 2, 6 \pmod{12}$, and $\dim(M_{k, 3}) = \lfloor k / 6 \rfloor + 1$ if $k \not\equiv 2, 6 \pmod{12}$. We have $M_{k, 3} = \mathbb{C} E_{k - 12 n, 3}^{*} E_{12 n, 3}^{*} \oplus M_{k, 3}^0$. Then,
\begin{equation*}
M_{k, 3}
 = E_{k - 12 n, 3}^{*} \left\{\mathbb{C} E_{12 n, 3}^{*} \oplus E_{12 (n-1), 3}^{*} M_{12, 3}^0 \oplus \cdots \oplus (M_{12, 3}^0)^n \right\} \oplus M_{k - 12 n, 3}^0 (M_{12, 3}^0)^n
\end{equation*}

In conclusion, we have the following proposition:
\begin{proposition}
Let $k \geqslant 4$ be an even integer. For every $f \in M_{k, 3}$, we have
\begin{equation}
\begin{split}
v_{i / \sqrt{3}}(f) \geqslant s_k &\quad(s_k=0, 1 \; \text{such that} \; 2 s_k \equiv k \pmod{4}),\\
v_{\rho_3}(f) \geqslant t_k &\quad(t_k=0, 1, 2, 3, 4, 5 \; \text{such that} \; - 2 t_k \equiv k \pmod{12}).
\end{split}
\end{equation}
In particular, if $f$ is a constant multiple of $E_{k, 3}^{*}$, then the equalities hold. \label{prop-bd_ord_3}
\end{proposition}\quad

\begin{center}
{\large Acknowledgement.}
\end{center}
The authors would like to thank Professor Eiichi Bannai for suggesting these problems as a master's course project.

\quad \\
{\it E-mail address}:

miezaki@math.kyushu-u.ac.jp (Tsuyoshi Miezaki),

nozaki@math.kyushu-u.ac.jp (Hiroshi Nozaki),

j.shigezumi@math.kyushu-u.ac.jp (Junichi Shigezumi).
\end{document}